\documentclass[12pt]{amsart}

\usepackage{xy, graphicx, color, hyperref,array, mathtools}
\usepackage{hyperref, pdftexcmds}

\usepackage{amsmath}
\usepackage{amssymb}
\usepackage{amsfonts,bm}
 
\usepackage{mathrsfs}

\usepackage{tikz,verbatim}
 
 \usepackage{dsfont}

\makeatletter
\def\subsection{\@startsection{subsection}{3}%
  \z@{.9\linespacing\@plus.7\linespacing}{.1\linespacing}%
  {\normalfont\bfseries}}
\makeatother 
 
 \xyoption{all} 
 
\title[Divisibility of Character Values]{The Divisibility of $\mathrm{GL}(n, q)$ character values}
 
 \author{Varun Shah}
 \author{Steven Spallone}

 \newtheorem{thm}{Theorem}[section]
\newtheorem{c.intro}[thm]{Corollary}
\newtheorem{lemma}[thm]{Lemma}
\newtheorem{prop}[thm]{Proposition}
\newtheorem{cor}{Corollary}[thm]

\theoremstyle{definition}
\newtheorem{remark}[thm]{Remark}
 
\newtheorem{defn}[thm]{Definition} 
 
\newcommand{\nc}{\newcommand}

\nc{\theorem}{\thm}

 \nc{\corollary}{\cor}
\nc{\mc}{\mathcal}
\nc{\mb}{\mathbb}
\nc{\mf}{\mathfrak}
\nc{\ul}{\underline}
\nc{\ol}{\overline}
\nc{\N}{\mb N}
\nc{\R}{\mb R}
\nc{\Z}{\mb Z}
\nc{\Q}{\mb Q}
\nc{\C}{\mb C}
\nc{\ms}{\mathscr }

\nc{\dmo}{\DeclareMathOperator}

\nc{\mat}[4]{
    \begin{pmatrix}
      #1 & #2 \\
      #3 & #4
    \end{pmatrix}
}

 \DeclareFontFamily{U}{cbgreek}{}
\DeclareFontShape{U}{cbgreek}{m}{n}{
        <-6>    grmn0500
        <6-7>   grmn0600
        <7-8>   grmn0700
        <8-9>   grmn0800
        <9-10>  grmn0900
        <10-12> grmn1000
        <12-17> grmn1200
        <17->   grmn1728
      }{}
 \DeclareRobustCommand{\qoppa}{%
  \text{\usefont{U}{cbgreek}{\normalorbold}{n}\symbol{19}}%
}

\makeatletter
\newcommand{\normalorbold}{%
  \ifnum\pdf@strcmp{\math@version}{bold}=\z@ bx\else m\fi
}
\makeatother

\dmo{\Ker}{Ker} \dmo{\val}{val} \dmo{\ord}{ord}

\dmo{\I}{I}
\dmo{\II}{II}
\dmo{\odd}{odd}
\dmo{\sgn}{sgn}
\nc{\beq}{\begin{equation*}}
\nc{\eeq}{\end{equation*}}
\nc{\half}{\frac{1}{2}}
\dmo{\lcm}{lcm}

\dmo{\Mod}{mod}

\dmo{\core}{core}
\dmo{\res}{res}
\dmo{\lin}{lin}
\dmo{\quo}{quo}
\dmo{\Sp}{Sp}

\dmo{\GL}{\mathrm{GL}}
 
\nc{\la}{\lambda}
  \nc{\eps}{\varepsilon}
  
 \nc{\lip}{\langle}
 \nc{\rip}{\rangle}
\nc{\gm}{\gamma}

\dmo{\Perm}{Perm}
\dmo{\Res}{Res}
\dmo{\Ind}{Ind}
\dmo{\tr}{tr}
\dmo{\Sym}{Sym}
\dmo{\reg}{reg}
\dmo{\End}{End}
\dmo{\even}{even}
\dmo{\Prop}{Prop}
\dmo{\irr}{Irr}
\dmo{\F}{\mathbf{F}}
 
\dmo{\supp}{supp}
\dmo{\barZ}{\overline{\Z}}
\dmo{\stab}{stab}
\dmo{\irrep}{\mathsf{Irr}}

\address{Department of Mathematics, University of Washington, Box 354350, Seattle, WA, 91895, USA}
\email{varunsh@uw.edu}
\address{Indian Institute of Science Education and Research, Pune-411008,Maharashtra,India}
\email{sspallone@iiserpune.ac.in}
\date{\today}


\begin{document}
	
	\maketitle
        \begin{center}
        \end{center}

 \begin{abstract}
     Let $q$ be a prime power, and $d$ a positive integer.  We study the proportion of irreducible characters of $\GL(n,q)$ whose values evaluated on a fixed matrix $g$ are divisible by $d$.  As $n$ approaches infinity, this proportion tends to $1$ when $q$ is coprime to $d$.  When $q$ and $d$ are not coprime, and $g=1$, this proportion is bounded above by $\frac{1}{q}$.  \end{abstract}

 \tableofcontents

 \section{Introduction}

 Let $q$ be a prime power, and $d$ a positive integer, and write $G_n=\GL(n,q)$. Fix a matrix $g \in G_{n_0}$, for $n_0 \leq n$. When $\chi$ is a character of $G_n$, we understand $\chi(g)$ to be $\chi$ evaluated on the block matrix 
 \beq \mat{g}{0}{0}{I_{n-n_0}} \in G_n. \eeq 
 In this paper we study the proportion of irreducible characters of 
 $G_n$ whose values at $g$ are divisible by $d$. 
(Here ``$d$ divides $\chi(g)$'' simply means that the algebraic number $\chi(g)/d$ is an algebraic integer.)
 Write $\irr(G_n)$ for the set of irreducible characters of $G_n$. Here is our main theorem.

 \begin{thm}\label{nuclear}
    If $d$ and $q$ are coprime,  then
    \beq 
    \lim \limits_{n \to \infty} \frac{\#\{\chi \in \irr(G_n) : d \mid \chi(g) \}}{|\irr(G_n)|} = 1.
    \eeq
 \end{thm}

Analogous results have been established for symmetric groups  in \cite{ganguly}, and for many simple Lie groups in \cite{Shah.Lie}. 

From \cite{Gallagher}, the proportion of \emph{zero} entries in the character table of $G_n$ approaches $100 \%$ as $n \to \infty$. Please note the different nature of the statistic: we are looking at divisibility of entries of a fixed column of the character table. (See also \cite{miller} for simple finite groups of Lie type, and \cite{pel3} for the symmetric groups.)

  On the other hand for $g$ the identity matrix, we have:
 
 \begin{thm}\label{voltas}
    If $d$ and $q$ are not coprime, then
    \beq
        \frac{\#\{\chi \in \irr(G_n): d \mid \deg \chi \}}{|\irr(G_n)|}  \leq \frac{1}{q}.
    \eeq
 \end{thm}

Our idea for Theorem \ref{nuclear} is as follows. Characters $\chi_{\bm \mu}$ of $G_n$ are indexed by partition-valued functions $\bm \mu$ of irreducible polynomials over ${\bf F}_q$. The degree of $\chi_{\bm \mu}$ is a multiple of a product of certain  $q$-analogues of the   hooklength formula for degrees of characters of the symmetric group  indexed by the values of $\bm \mu$. The problem for $g=I$ then reduces to the asymptotic results of \cite{ganguly}. Finally we  use a well-known relation (Lemma  \ref{was.previous} below) between character values and character degrees.

 In Section \ref{notate.sec} we establish notations and review \cite{ganguly}, and the degree formula for $G_n$. In Section  \ref{valuations.sec}  we use Lemma \ref{was.previous} to give a valuative criteria for the divisibility of $\chi(g)$ by $d$. In Section \ref{divisibility.powers.two}  we prove Theorem \ref{nuclear} when $d$ is a power of $2$. The case of odd prime divisors is  more intricate, involving the core and quotient theory of partitions, but the same overall approach succeeds.  We carry this out in Section \ref{odd.div.section}.   Theorem \ref{voltas} is an easy exercise contained  in Section \ref{def.char.sec}.
 

\bigskip

\textbf{Acknowledgements.} 
The authors would like to thank Jeff Adler for helpful discussions, and Sandeep Varma for pointing us to Lemma \ref{was.previous}.

\section{Notations and Preliminaries} \label{notate.sec}
Write $\ul \N$ for the set of nonnegative integers.
 \subsection{Partitions} \label{part.section}
    
\begin{defn}
    A \emph{partition} $\lambda = \lambda_1 \geq \cdots \geq \lambda_k > 0$ is a weakly decreasing sequence of positive integers. 
\end{defn}
Write $|\lambda| = \sum_i \lambda_i$ for the size of $\la$, and $\ell(\la)=k$ for the length of $\la$.  We write $\Lambda$ for the set of partitions, and $\Lambda_n$ for those of size $n$. As usual write $p(n)=|\Lambda_n|$, the number of partitions of $n$. Write $\mc H(\la)$ for the set of hooks of $\la$.    

 We will need the theory of $t$-cores and $t$-quotients of partitions, for a positive integer $t$. Their definitions and properties can be found in \cite{Robinson}, for instance. 
 Here we simply state the key properties for later reference. A hook $h$ of $\la$ is called a \emph{$t$-hook} when $|h|$ is a multiple of $t$. Write $\mc H_t(\la)$ for the set of $t$-hooks. When $\la$ does not have any $t$-hooks, it is called a \emph{$t$-core partition}. Write $\mc C^t$ for the set of $t$-core partitions, and $\mc C^t_n$ for those of size $n$.

The \emph{$t$-core} of a partition $\la$, written $\core_t(\la)$, is a certain $t$-core partition contained in $\la$. The \emph{$t$-quotient} of $\la$, written $\quo_t \la=(\la^0, \ldots, \la^{t-1})$, is a certain $t$-tuple of partitions. 
There is a bijection between $\mc H_t(\la)$ and the disjoint union
\beq
\mc H(\la^0) \amalg \cdots \amalg \mc H(\la^{t-1}),
\eeq
under which the sizes of the hooks are divided by $t$. Thus 
\beq
\sum_{i=0}^{t-1} |\la^i|=|\mc H_t(\la)|.
\eeq
Moreover, the size of $\la$ can be recovered by
\begin{equation} \label{zepto}
|\la|=|\core_t \la|+ t \sum_{i=0}^{t-1} |\la^i|.
 \end{equation}
The map
 \beq
 \Lambda \overset{\sim}{\to} \underbrace{\Lambda \times \cdots \times \Lambda}_{t \text{ times}} \times \: \mc C^t
 \eeq
 given by $\la \mapsto (\la^0, \ldots, \la^{t-1},\core_t(\la))$ is bijective.

\subsection{Divisibility of \texorpdfstring{$S_n$}{Sn} Characters} 

We need an analogue of Theorem \ref{nuclear} for symmetric group characters in \cite{ganguly}. Irreducible characters of $S_n$ correspond to partitions of size $n$. Denote by $\mf f_\lambda$ the degree of the character corresponding to $\lambda \vdash n$; it is given by 
\begin{equation} \label{hl.formula}
\mf f_\la=\frac{n!}{\prod_{h \in \mc H(\la)} |h|}.
\end{equation}

\begin{prop}\label{band} \cite[Theorem A]{ganguly}
     For every prime number $\ell$ and   $r>0$,
     \beq
        \lim\limits_{n \to \infty} \frac{\#\{\lambda \vdash n \mid v_\ell(\mf f_\lambda) < r +  \log_\ell n\}}{p(n)} = 0.
     \eeq
\end{prop}

We shall also use the shorthand
\beq
f_{\quo_t \la}= \prod_{i=0}^{t-1} f_{\la^i}.
\eeq

\subsection{Character Degrees for $G_n$} \label{char.deg.sec}

In this section we recall the character degree formulas for $G_n$, in terms of Green's parameterization, and give a convenient factorization.
  
\subsubsection{Degree Formula}\label{hanger}

 Let $n$ be a positive integer and $q$ a prime power. 
 Write  $\Phi$ for the set of monic irreducible polynomials in $\F_q[ x]$ which are not equal to $f(x)=  x$. For $f \in \Phi$, put $d(f)=\deg f$.   Let $\mf X_{n}$ be the set of functions ${\bm \mu}\colon \Phi \to \Lambda$ with
 \beq
    \sum\limits_{f \in \Phi} |\bm \mu(f)|d(f) = n.
 \eeq
 The irreducible characters of $G_n$ correspond to elements of $\mf X_{n}$ and were first explicitly computed   in \cite{green}. Let $\chi_{\bm \mu}$ be the character corresponding to ${\bm \mu} \in \mf X_{n}$, and put $d_{\bm \mu} =\deg \chi_{{\bm \mu}}$.
 
 For a partition $\la$, write $\alpha(\lambda) = \sum_i(i-1)\lambda_i$. According to \cite[Section 8.7]{springer},  
 \begin{equation} \label{remote}
    d_{\bm \mu} = \psi_n(q)\prod\limits_{f \in \Phi}H\left(\bm \mu(f), q^{d(f)}\right),
 \end{equation}
where 
\beq
\psi_n(q) = (q^n-1)(q^{n-1}-1)\cdots(q-1),
\eeq
and for a partition $\lambda$, 
\beq 
H(\lambda, x) = x^{\alpha(\lambda)}\prod\limits_{h \in \mc H(\lambda)}\left(x^{|h|} - 1\right)^{-1}.
\eeq
 Our next step is to express $d_{\bm \mu}$ in terms of degrees of unipotent characters.

 \subsubsection{Unipotent Characters}

 For each partition $\lambda \vdash n$, we  define $\bm \mu_{\la} \in \mf X_n$ 
by ${\bm \mu_{\la}}(x-1)=\lambda$, and $\bm \mu_{\la}(f)=\emptyset$ for $f(x) \neq x-1$. Write $\chi_\la$ for $\chi_{\bm \mu_{\la}}$; we call such characters \emph{unipotent}. 

\begin{remark} Unipotent characters are the irreducible constituents of the permutation representation of $G_n$ on $G_n/B_n$, where $B_n$ is the subgroup of upper triangular members.
\end{remark}

Let $d_\la=\deg(\chi_\la)$. From  \eqref{remote}, this is
\begin{equation} \label{cake}
    d_\lambda=d_\la(q) =  q^{\alpha(\la)}\frac{\prod \limits_{i=1}^{n}q^i - 1}{\prod\limits_{h \in \mc H(\lambda)}q^{|h|} - 1}.
\end{equation}

\begin{remark}
This may be regarded as a $q$-analogue of \eqref{hl.formula}, since $\mf f_\la$ is the limit  as $q \to 1$.
\end{remark}

\subsubsection{Factoring the Degree}

In what follows, it will be convenient to view $d_\la(q)$ in \eqref{cake} as a function of $q$. For $\bm \mu \in \mf X_{n}$, we may write
\beq
\begin{split}
H(\bm \mu(f),q^{d(f)}) &= q^{\alpha(\bm \mu(f))d(f)} \prod_{h \in \mc H(\bm \mu(f))} (q^{d(f)|h|}-1)^{-1} \\
			&= \frac{d_{\bm \mu(f)}(q^{d(f)})}{\psi_{|\bm \mu(f)|}(q^{d(f)})}. \\
			\end{split}
			\eeq
From   \eqref{remote}, we may write
\beq
\begin{split}
d_{\bm \mu}&=\frac{\psi_n(q)}{\prod_f \psi_{|\bm \mu(f)|}(q^{d(f)})} \cdot \prod_f d_{\bm \mu(f)}(q^{d(f)}) \\
 &= a_{\bm \mu} \cdot b_{\bm \mu}, \\
\end{split}
\eeq
  where
\begin{equation} \label{switch}
a_{\bm \mu} = \frac{\prod\limits_{i = 1}^n q^i-1}{\prod\limits_{f \in \Phi}\prod\limits_{i = 1}^{|\bm \mu(f)|}q^{d(f)i}-1}
\end{equation}
and
\beq
b_{\bm \mu}= \prod\limits_{f\in \Phi}d_{\bm \mu(f)}(q^{d(f)}).
\eeq

Thus $a_{\bm \mu}$ only depends on the sizes of the partitions $\bm \mu(f)$, and $b_{\bm \mu}$ is a product of degrees of unipotent characters evaluated at powers of $q$.

\begin{remark} This factorization is motivated by Lusztig's Jordan Decomposition of characters. In  \cite[Corollary 2.6.6]{Geck},  our $b_{\bm \mu}$ is the degree of the ``unipotent factor'' of $\chi_{\bm \mu}$, and our $a_{\bm \mu}$ is the ``index factor''. 
\end{remark}

 \section{Divisibility of Character Values} \label{valuations.sec} 
\subsection{Valuations}

Being a sum of roots of unity, a character of a finite group takes values in the ring of integers $\mc O$ of a number field $K$. 
We recall here the theory of  valuations as described in \cite[Chapter 11]{Janusz}.
Given a prime ideal $\mf p$ of $\mc O$, and $x \in \mc O$, write $v_{\mf p}(x)$ for the highest exponent $r$ so that $x \in \mf p^r$.  
This extends to a valuation of $K$ in the usual way, via $v_{\mf p}(a/b)=v_{\mf p}(a)-v_{\mf p}(b)$ for $a,b \in \mc O$ with $b \neq 0$. For example, if $K=\Q$, and $\mf p=\ell\Z$, then $v_{\mf p}$ is the usual $\ell$-adic valuation $v_\ell$, meaning that for $n$ an integer, $v_\ell(n)$  is the highest exponent $r$ so that $\ell^r$   divides $n$. 

Let $\ell$ be the prime in $\Z$ lying under $\mf p$, meaning $\mf p \cap \Z=\ell \Z$. 
Then there is a positive integer $e=e(\mf p)$ (the ``ramification index") so that for all $x \in \Q$, we have $v_{\mf p}(x)=e \cdot v_\ell(x)$, by \cite[Theorem 3.1]{Janusz}.

According to \cite[Proposition 48, page 720]{DF}, $x \in K$ is integral iff   $v_{\mf p}(x) \geq 0$ for all $\mf p$. Therefore:

\begin{lemma} \label{fixit} Let $\alpha,\delta \in \mc O$, with $\delta \neq 0$. If for all $\mf p \ni \delta$, we have $v_{\mf p}(\alpha) \geq v_{\mf p}(\delta)$, then $\frac{\alpha}{\delta} \in \mc O$.
\end{lemma}

 \subsection{Relating Character Degrees to Character Values}

When an irreducible character has a highly divisible degree, its other values will also be highly divisible, by the following.

\begin{lemma}  \cite[Exercise 6.9]{serre.linear},  \label{was.previous}
    Let $G$ be a finite group, $g \in G$ and $\chi$ an irreducible   character of $G$. Then 
    \beq
        \frac{\chi(g)}{\deg \chi}[G: Z_G(g)]
    \eeq
    is an (algebraic) integer.
\end{lemma}

In our situation, the centralizer $Z_n(g)$ of $g$ in $G_n$ contains $G_{n-n_0}$ as a subgroup via block matrices.
The lemma asserts that 
    \beq
    \frac{\chi_{\bm \mu}(g)}{d_{\bm \mu}} [G_n: Z_n(g)]
    \eeq
    is an integer, hence so too is
     \beq
    \frac{\chi_{\bm \mu}(g)}{d_{\bm \mu}} [G_n: G_{n-n_0}].
    \eeq 
    We deduce:

\begin{lemma}
 For all $\bm \mu \in \mf X_{n}$ and $g \in G_{n_0} \leq G_n$, the quantity
\beq
\frac{\chi_{\bm \mu}(g)  \cdot q^{\binom{n}{2} - \binom{n-n_0}{2}}   \prod\limits_{i=0}^{n_0-1}    q^{n-i}-1   }{d_{\bm \mu}} 
\eeq
is an integer.
\end{lemma} 

\begin{cor}\label{helmet} 
Let $d$ be coprime to $q$. If for each prime divisor $\ell$ of $d$, we have 
\beq
v_\ell(d_{\bm \mu}) -v_\ell \left(   \prod\limits_{i=0}^{n_0-1}    q^{n-i}-1 \right) \geq v_\ell(d),
\eeq
then $d$ divides $\chi_{\bm \mu}(g)$.
\end{cor}
 
\begin{proof}
 We will apply Lemma \ref{fixit}.
Let 
\beq
\beta=\frac{\chi_{\bm \mu}(g)  \cdot q^{\binom{n}{2} - \binom{n-n_0}{2}}   \prod\limits_{i=0}^{n_0-1}    q^{n-i}-1   }{d_{\bm \mu}} \in \mc O
\eeq
and 
\beq
c=\frac{d_{\bm \mu}}{d q^{\binom{n}{2} - \binom{n-n_0}{2}}   \prod\limits_{i=0}^{n_0-1}    q^{n-i}-1   } \in \Q,
\eeq
so that $\frac{\chi_{\bm \mu}(g)}{d}=\beta \cdot c$. For a prime ideal $\mf p$ of $\mc O$, we have
\beq
v_{\mf p} \left( \frac{\chi_{\bm \mu}(g)}{d} \right) =v_{\mf p}(\beta)+v_{\mf p}(c) \geq v_{\mf p}(c).
\eeq

If $d \in \mf p$ and $\ell$ lies under $\mf p$, then $\ell$ divides $d$. Since $\gcd(d,q)=1$, we cannot have $\ell | q$, hence $v_\ell(q)=0$. Now the hypothesis implies that $v_{\ell}(c)  \geq 0$. So $v_{\mf p}(c)=e  \cdot v_{\ell}(c) \geq 0$, where  $e=e(\mf p)>0$ is the ramification index. Thus $v_\mf p(\chi_{\bm \mu}(g)) \geq v_{\mf p}(d)$ for all $\mf p \ni d$, and the conclusion follows from Lemma \ref{fixit}. 
 
\end{proof}

  \section{Divisibility by Powers of $2$}\label{divisibility.powers.two}

\subsection{A Lower Bound for the $2$-Divisibility of the Degree}
In this section, let $q$ be an odd prime, and write  $v$ for the $2$-adic valuation $v_2$. We will repeatedly use the following  elementary number theory fact.
\begin{lemma} (Lifting-the-exponent or LTE, \cite{Cuellar}) \label{angles}
    For an odd integer $a$, 
    \beq
        v(a^n - 1) = \begin{cases}
v(a-1) + v(a+1) + v(n) - 1, & n \text{ even} \\	
v(a-1), & n \text{ odd}
\end{cases}.
    \eeq
\end{lemma}

\begin{prop} \label{socks}
    We have
    \beq
      v(d_\lambda) \geq v(\mf f_\lambda).
    \eeq
\end{prop}

\begin{proof}

By LTE,
\beq
 v\left(\prod \limits_{i=1}^{n}q^i - 1\right)=n v(q-1)+ \left\lfloor \frac{n}{2} \right\rfloor v \left(\frac{q+1}{2} \right) + \sum_{i \even} v(i).
 \eeq
Hence by \eqref{cake} and another application of LTE:
\beq
\begin{split}
    v(d_\lambda) &= v\left(\prod \limits_{i=1}^{n}q^i - 1\right) - v\left(\prod\limits_{h \in \mc H(\lambda)}q^{|h|} - 1\right) \\
    &= \left(\left\lfloor\frac{n}{2}\right\rfloor - \#\{h \in \mc H(\lambda) \colon |h| \text{ even}\}\right)v\left(\frac{q+1}{2}\right) + v(\mf f_\lambda)  \\
    & \geq v(\mf f_\la) \text{ by }  \eqref{zepto}. \\
\end{split}
\eeq \end{proof}

By Proposition \ref{socks}, we see that
\begin{equation} \label{b-est}
v(b_{\bm \mu}) \geq \sum\limits_{f \in \Phi}v\left(\mathfrak f_{\bm \mu(f)}\right).
\end{equation}

\begin{lemma} \label{adidas}
 We have
    \beq
        v(a_{\bm \mu}) \geq v\left(\frac{n!}{\prod\limits_{f\in \Phi} |\bm \mu(f)|!}\right).
    \eeq
\end{lemma}

\begin{proof}
Applying $v$ to \eqref{switch} and using LTE,  
\begin{multline*}
v(a_{\bm \mu}) = \left[nv(q-1) + \left\lfloor\frac{n}{2}\right\rfloor v\left(\frac{q+1}{2}\right) + v(n!) \right] \\
- \sum\limits_{f\colon d(f) \text{ even}} \left[|\bm \mu(f)|v(q-1) + |\bm \mu(f)|v\left(\frac{q+1}{2}\right) + |\bm \mu(f)|v(d(f)) + v(|\bm \mu(f)|!)\right] \\
- \sum\limits_{f\colon d(f) \text{ odd}} \left[|\bm \mu(f)|v(q-1) + \left\lfloor\frac{|\bm \mu(f)|}{2}\right\rfloor v\left(\frac{q+1}{2}\right) + |\bm \mu(f)|v(d(f)) + v(|\bm \mu(f)|!)\right].
\end{multline*}

Since $n = \sum\limits_{f \in \Phi} d(f)|\bm \mu(f)|$, we can write
\beq
\begin{split}
    n &= \sum\limits_{f\colon d(f) \text{ even}} d(f)|\bm \mu(f)| + \sum\limits_{f\colon d(f) \text{ odd}} d(f)|\bm \mu(f)| \\
    &\geq \sum\limits_{f\colon d(f) \text{ even}} 2|\bm \mu(f)| + \sum\limits_{f\colon d(f) \text{ odd}} |\bm \mu(f)|.
\end{split}
\eeq

This means that the coefficient of $v\left(\frac{q+1}{2}\right)$, 
\beq
\left\lfloor\frac{n}{2}\right\rfloor - \sum\limits_{f\colon d(f) \text{ even}} |\bm \mu(f)| - \sum\limits_{f\colon d(f) \text{ odd}} \left\lfloor\frac{|\bm \mu(f)|}{2}\right\rfloor,
\eeq
is nonnegative. So

\beq 
v(a_{\bm \mu}) \geq \left(n - \sum\limits_{f \in \Phi}|\bm \mu(f)|\right)v(q-1) - \sum\limits_{f \in \Phi}|\bm \mu(f)|v(d(f)) + v\left(\frac{n!}{\prod\limits_{f\in \Phi} |\bm \mu(f)|!}\right).
\eeq
Similarly, since $n = \sum\limits_{f \in \Phi  } d(f)|\bm \mu(f)|$, we have 
\beq
n - \sum\limits_{f\in \Phi} |\bm \mu(f)| = \sum\limits_{f\in \Phi} (d(f)-1)|\bm \mu(f)|.
\eeq
This gives
\beq
v(a_{\bm \mu}) \geq \sum\limits_{f \in \Phi} [d(f)-1-v(d(f))]|\bm \mu(f)| + v\left(\frac{n!}{\prod\limits_{f \in \Phi} |\bm \mu(f)|!}\right),
\eeq
as $v(q-1) \geq 1$. Since for a positive integer $x$, 
\beq
    v(x) \leq \log_2(x) \leq x-1,
\eeq
the sum $\sum_f [d(f)-1-v(d(f))]|\bm \mu(f)|$ is nonnegative. The   conclusion follows.
 
\end{proof}

From Lemma \ref{adidas} and the inequality \eqref{b-est}, we deduce the following.

\begin{thm}
For $\bm \mu \in \mf X_{n}$,
\beq
v(d_{\bm \mu}) \geq v\left(\frac{n!}{\prod\limits_{f \in \Phi} |\bm \mu(f)|!}\right) + \sum\limits_{f \in \Phi}v\left(\mathfrak f_{\bm \mu(f)}\right).
\eeq
\end{thm}

\subsection{Asymptotics for the $2$-Divisibility of Degrees}

We continue to assume that $q$ is odd. Given integers $0<k \leq n$, let
\beq
    (n)_k = \frac{n!}{(n-k)!};
\eeq
this is called the \emph{falling factorial}.

The following proposition is somewhat stronger than the statement that ``$100\%$ of the irreducible characters of $G_n$ have degree divisible by any fixed power of $2$''.

\begin{prop}\label{bottle2}
    Let $k$ be a positive integer. Then for all $\eps > 0$, there is an $N = N(\eps)$ such that for all $n \geq N$,
    \beq
        \frac{\#\{\bm \mu \in \mf X_{n}: v(d_{\bm \mu}) <  v((n)_k)\}}{\#\mf X_{n}} < \eps.
    \eeq
\end{prop}

Let $\mathscr F_n(\Phi)$ be the set of functions $F$ from $\Phi$ to the set of whole numbers, with $\sum_{f \in \Phi} d(f)F(f)=n$. 
The \emph{support} of $F \in \ms F_n(\Phi)$, written $\supp(F)$, is the set of $f \in \Phi$ with $F(f) \neq 0$.
We have a surjection $\mf X_{n} \twoheadrightarrow \ms F_n(\Phi)$ written
 $\bm \mu \mapsto F_{\bm \mu}$, given by $F_{\bm \mu}(f)=|\bm \mu(f)|$.

 \begin{lemma} \label{blanket}
    Let $k$ be a positive integer. Then for all $\eps > 0$, there is an $N = N(\eps)$ such that for all $n \geq N$ and for all $F \in \ms F_n(\Phi)$ we have

    \beq
        \frac{\#\{\bm \mu: \bm \mu \mapsto F, \text{ and } v(d_{\bm \mu}) <  v((n)_k)\}}{\#\{\bm \mu: \bm \mu \mapsto F\}} < \eps.
    \eeq
\end{lemma}

\begin{proof}
    
Fix $\eps > 0$. By Proposition \ref{band} there is an $M = M(\eps)$ such that for all $m \geq M$, 
\beq 
    \frac{\#\{\lambda \vdash m : v(\mf f_\lambda) <  k+ \log_2(kq^{k+1}) + \log_2m \}}{p(m)} < \eps.
\eeq

Let $N = kq^{k+1}M$ and $F \in \ms F_n(\Phi)$ for $n \geq N$. Suppose $\max\limits_{f \in \supp F}d(f) \leq k$. As there are $q^{k+1}$ polynomials of degree at most $k$ in $\F_q[x]$, we have $|\supp F| \leq q^{k+1}$. Thus
\beq
    \begin{split}
        \frac{1}{|\supp F|}\sum\limits_{f \in \Phi} F(f) &\geq \frac{\sum\limits_{f \in \Phi} F(f)}{q^{k+1}} \\
        &\geq \frac{\sum\limits_{f \in \Phi} d(f)F(f)}{kq^{k+1}} \\
        &\geq \frac{n}{kq^{k+1}} \\
        &\geq M.
    \end{split}
\eeq

So, there is an $f_0 \in \Phi$ with $F(f_0) \geq M$, say $m=F(f_0)$. For any $\bm \mu \mapsto F$ we have $v(d_{\bm \mu}) \geq v(\mf f_{\bm \mu(f_0)})$. Since the fibre over $F$ corresponds to products of partitions of each $F(f)$, we have
\beq
\begin{split}
    \frac{\#\{\bm \mu: \bm \mu \mapsto F: v(d_{\bm \mu})< v((n)_k)}{\#\{\bm \mu: \mb \mu \mapsto F\}} &= \frac{ \#\{ \bm \mu \in \prod_{f \in \supp F} \Lambda_{F(f)} \mid v(d_{\bm \mu})  <  v((n)_k)\}}{ \prod_f p(F(f))} \\
    &\leq \frac{ \#\{ \bm \mu \in \prod_{f \in \supp F} \Lambda_{F(f)} \mid v(\mf f_{\bm \mu(f_0)})  <  v((n)_k)\}}{ \prod_f p(F(f))} \\
&= \frac{\#\{\lambda \vdash m : v(\mf f_\lambda) <  v((n)_k)\}}{p(m)} \\
&\leq \frac{\#\{\lambda \vdash m : v(\mf f_\lambda) <  k+ \log_2(kq^{k+1})+ \log_2m }{p(m)}< \eps. \\
   \end{split} 
\eeq
On the other hand, if $\max\limits_{f \in \supp F}d(f) > k$, let $f_1 \in \supp F$ have the largest degree. Then  
\beq
\begin{split}
n &= \sum \limits_{f\in \Phi}d(f) F(f) \\
&=  d(f_1)F(f_1)+ \sum\limits_{f \neq f_1}d(f)F(f) \\
&> kF(f_1)+ \sum\limits_{f \neq f_1}F(f)  \\
&= (k-1)F(f_1)+\sum\limits_{f \in \Phi}F(f). \\
\end{split}
\eeq

Because $F(f_1) \geq 1$, we have $n - k  \geq \sum\limits_{f \in \Phi}F(f)$,
and therefore the quotient  $\dfrac{n!}{\prod\limits_{f\in \Phi}F(f)!}$ is a product of $(n)_k$ with certain multinomial coefficients.
For any $\bm \mu \mapsto F$, we have
\beq
v(d_{\bm \mu}) \geq v\left(\frac{n!}{\prod\limits_{f\in \Phi}F(f)!}\right) \geq v((n)_k). 
\eeq
So, in this case 
\beq
    \frac{\#\{\bm \mu: \bm \mu \mapsto F \text{ and } v(d_{\bm \mu}) <  v((n)_k)\}}{\# \{\bm \mu: \bm \mu \mapsto F\}} = 0. 
\eeq \end{proof}
The proposition now simply follows:
\begin{proof}[Proof of Proposition \ref{bottle2}]
    For any $\eps>0$, let $N=N(\eps)$ be as above.  Then by Lemma \ref{blanket},
    \beq
    \begin{split}
        \frac{\#\{\bm \mu \in \mf X_{n}: v(d_{\bm \mu}) <  v((n)_k)\}}{|\mf X_{n}|} &= \sum_{F } \frac{\#\{\bm \mu \mapsto F: v(d_{\bm \mu}) <  v((n)_k)\}}{|\mf X_{n}|} \\
         &< \eps\sum_{F }\frac{\#\{ \bm \mu:\bm \mu \mapsto F\}}{|\mf X_{n}|} \\
        &= \eps,
    \end{split}
    \eeq
    the sums being over $F \in \ms F_n(\Phi)$.
\end{proof}

\subsection{$2$-Divisibility of Character Values} \label{2.div.char.sec}

Now we prove Theorem \ref{nuclear} in the case where $d$ is a power of $2$.

\begin{prop} \label{bip} Let   $r$ be a positive integer. Then
\beq 
        \lim\limits_{n \to \infty} \frac{\#\{\bm \mu \in \mf X_{n}:  2^r \nmid   \chi_{\bm \mu}(g)\}}{|\mf X_{n}|} = 0.
    \eeq
\end{prop}


\begin{proof} 
     By LTE, we have
    \beq
        v\left(\prod\limits_{i=0}^{n_0-1} q^{n-i}-1\right) \leq n_0v\left(q^2-1\right) + v((n)_{n_0}).
    \eeq
    Corollary \ref{helmet} implies that if $2^r \nmid \chi_{\bm \mu}(g)$, then 
    \beq
        v(d_{\bm \mu}) < r + v\left(\prod\limits_{i=0}^{n_0-1} q^{n-i}-1\right)  \leq  r + n_0v\left(q^2-1\right) + v((n)_{n_0}).
    \eeq
    Let $r_0 = r + n_0 \: v(q^2-1)$. Then
    \beq
    v(d_{\bm \mu}) < r_0 + v((n)_{n_0}) \leq v((n)_{n_0+2r_0}),
  \eeq  
     whenever $2^r \nmid \chi_{\bm \mu}(g)$. So, by Proposition \ref{bottle2},
    \beq
        \lim\limits_{n \to \infty} \frac{\#\{\bm \mu \in \mf X_{n}:  2^r  \nmid \chi_{\bm \mu}(g)\}}{|\mf X_{n}|} = 0.
    \eeq
\end{proof}

\section{Odd Divisors} \label{odd.div.section}

In this section, let $\ell$ be an odd prime coprime to $q$. We now consider divisibility of character values by powers of $\ell$, so simply write $v$ for the $\ell$-adic valuation $v_\ell$.
 Write $t$ for the multiplicative order of $q$ mod $\ell$, and put $\tau =v(q^t-1) \geq 1$.
\subsection{Lifting the Exponent Again}

We need LTE for odd primes.

\begin{lemma}  \cite{Cuellar}
   If $n$ is a multiple of $t$, then
    \beq
        v(q^n - 1) = v(n/t)+ \tau.
        \eeq
\end{lemma}

\begin{lemma} \label{c.calcs}
Given positive integers $A,B$, 
\beq
v \left( \prod_{i=1}^A q^{Bi}-1 \right)=\left\lfloor Ah/t   \right\rfloor (v(B/h)+ \tau)+ v( \left\lfloor Ah/t  \right\rfloor!),
\eeq
where $h=\gcd(B,t)$.
\end{lemma}

\begin{proof} Let $m=\lcm(B,t)$. Then 
\begin{equation*}
\begin{split}
v \left( \prod_{i=1}^A q^{Bi}-1 \right) &= v \left( \prod_{i=1}^{ \left\lfloor Ah/t  \right\rfloor} q^{im}-1 \right) \\
&= \sum_{i=1}^{ \left\lfloor Ah/t  \right\rfloor} \tau + v \left(\frac{iB}{h} \right) \\
&= \left\lfloor Ah/t  \right\rfloor \tau+ v \left(\prod_i  \frac{iB}{h} \right)\\
&= \left\lfloor Ah/t  \right\rfloor(\tau+ v(B/h))+ v( \left\lfloor Ah/t  \right\rfloor!).\\
\end{split}
\end{equation*} \end{proof}

Note the   case
\begin{equation} \label{cn1}
v \left( \prod_{i=1}^n q^i-1 \right)   =\left \lfloor n/t \right \rfloor \tau+ v(\left \lfloor n/t \right \rfloor !). \\
\end{equation}

\subsection{The $\ell$-adic Valuation of the Degree}

We need an analogue of Proposition \ref{socks}.
\begin{prop}  \label{monsoonz} 
    For a partition $\la$,
    \beq
      v(d_\lambda) \geq \left \lfloor \frac{|\core_t(\la)|}{t} \right \rfloor + v(\mf f_{\quo_t(\la)}).
    \eeq
\end{prop}


\begin{proof}  
We have
\beq
 v \left( \prod_{h \in \mc H(\la)} q^{|h|}-1 \right) = \# \{ \text{$t$-hooks} \}\tau + \sum_{\text{$t$-hooks } h} v \left(\frac{|h|}{t} \right). 
 \eeq

Thus from \eqref{cake} and \eqref{cn1},
\beq
v (d_\la) = \left( \left\lfloor n/t \right\rfloor- \# \{\text{$t$-hooks}\} \right) \tau+v \left( \frac{\left\lfloor n/t \right\rfloor!}{\prod_{\text{$t$-hooks}}\frac{|h|}{t}} \right). 
\eeq
By the properties of core and quotient (see Section \ref{part.section}) and \eqref{hl.formula}, we deduce that
\beq
v (d_\la)=  \left\lfloor   \frac{|\core_t(\la)|}{t} \right \rfloor \tau + v \left( \frac{ \left\lfloor n/t \right\rfloor!}{|\la^0|! \cdots |\la^{t-1}|!} \mf f_{\la^0} \cdots \mf f_{\la^{t-1}} \right).
 \eeq \end{proof}

For $f \in \Phi$, put $h_f=\gcd(d(f),t)$. Recall the integer $a_{\bm \mu}$ defined by \eqref{switch}. Here is our replacement for Lemma \ref{adidas}. 

\begin{lemma} For any ${\bm \mu} \in \mf X_{n}$,
\beq
v(a_{\bm \mu}) \geq v \left( \frac{\left\lfloor n/t \right\rfloor!}{\prod_{f \in \Phi} \left \lfloor \frac{|\bm \mu(f)| h_f}{t} \right \rfloor!} \right).
\eeq
\end{lemma}

\begin{proof} 
By Lemma \ref{c.calcs},
\beq
v(a_{\bm \mu}) =\left\lfloor n/t \right\rfloor\tau+v(\left\lfloor n/t \right\rfloor!)-\sum_f \left \lfloor \frac{|\bm \mu(f)|h_f}{t} \right \rfloor \left(\tau+ v \left(\frac{d(f)}{h_f} \right) \right)+ v \left( \left \lfloor \frac{|\bm \mu(f)|h_f}{t} \right \rfloor \right),
\eeq
so it is enough to show that
\begin{equation} \label{goal.here}
\left\lfloor n/t \right \rfloor \tau \geq \sum_f \left \lfloor \frac{|\bm \mu(f)|h_f}{t} \right \rfloor \left(\tau+ v \left(\frac{d(f)}{h_f} \right) \right).
 \end{equation}
Now since
\beq
n=\sum_f d(f) |\bm \mu(f)|,
\eeq
  we can write
\beq
\frac{n}{t}=\sum_f \frac{|\bm \mu(f)|}{t} (d(f)-h_f)+\frac{|\bm \mu(f)|h_f}{t},
\eeq
and therefore
\begin{equation} \label{staplerz} 
 \left\lfloor n/t \right \rfloor  \tau \geq \sum_f \left \lfloor \frac{|\bm \mu(f)|}{t} (d(f)-h_f) \right \rfloor+  \left \lfloor \frac{|\bm \mu(f)|h_f}{t} \right \rfloor \tau.
\end{equation} 
Generally, for an integer $x$, we have $x \geq v(x)+1$. Setting $x=\frac{d(f)}{h_f}$ we deduce
\beq
\frac{|\bm \mu(f)|}{t} (d(f)-h_f) \geq \frac{|\bm \mu(f)|h_f}{t} v \left(\frac{d(f)}{h_f} \right),
\eeq
and so \eqref{goal.here} follows from  \eqref{staplerz}.
\end{proof} 
 
 From Proposition \ref{monsoonz},
\beq
v(b_{\bm \mu}) \geq \sum_{f \in \Phi} \left( \left \lfloor \frac{|\core_t(\bm \mu(f))|}{t} \right \rfloor + v(\mf f_{\quo_t \bm \mu(f)}) \right),
\eeq
and so altogether,
\begin{equation} \label{needed}
v(d_{\bm \mu}) \geq v \left(   \frac{\left\lfloor n/t \right \rfloor!}{\prod_{f \in \Phi} \left \lfloor \frac{|\bm \mu(f)| h_f}{t} \right \rfloor!}   \right)+
 \sum_{f \in \Phi} \left( \left \lfloor \frac{|\core_t(\bm \mu(f))|}{t} \right \rfloor + v(\mf f_{\quo_t \bm \mu(f)}) \right).
 \end{equation}

Consider the map $\qoppa: \Lambda \to \ul \N^{t+1}$ (``qoppa'', the Greek $q$) given by
\beq
\qoppa(\la)=\left( t|\la^0|, \cdots, t|\la^{t-1}|,  |\core_t(\la) | \right).
\eeq
We put $|(n_0, \ldots, n_t)|=n_0+\ldots +n_{t-1}+n_t$, so that $|\qoppa(\la)|=|\la|$. Let $\ms F(\Phi)^{t+1}$ be the set of functions  $\mathsf F:\Phi \to \ul \N^{t+1}$.
Let us write $\mathsf F(f)=(\mathsf F_0(f), \mathsf F_1(f), \ldots, \mathsf F_t(f))$, and $||\mathsf F||=\sum_f d(f) |\mathsf F(f)|$.

Let $\mf I_n$ be the image of the map $\mf X_n \to \ms F(\Phi)^{t+1}$  defined by $\bm \mu \mapsto  \qoppa \circ \bm \mu$. Then $\mf I_n$ is the set of all $\mathsf F$ with the following three properties:
\begin{itemize}
\item For all $f \in \Phi$, there is a $t$-core partition of size $\mathsf F_t(f)$.
\item For all $f \in \Phi$ and integers $0 \leq i < t$, the integer $t$ divides $\mathsf F_i(f)$.
\item $||\mathsf F||=n$.
\end{itemize}

 It is convenient to extend the falling factorial to positive numbers. For $0<r<s$, we set $
(s)_r=\frac{\lfloor s \rfloor!}{\lfloor s-r \rfloor!}$.
Legendre's formula implies that
\beq
v((s)_r) \leq \lceil r \rceil+\log_\ell(s).
\eeq

 \begin{lemma} \label{blanket2}
For all $u,\eps > 0$, there is an $N = N(\eps)$ such that for all $n \geq N$ and for all $\mathsf F \in \mf I_n$ we have

    \beq
        \frac{\#\{\bm \mu: \bm \mu \mapsto \mathsf F, \text{ and } v(d_{\bm \mu}) < v((n/t)_u)\}}{\#\{\bm \mu: \bm \mu \mapsto \mathsf F\}} < \eps.
    \eeq
\end{lemma}

\begin{proof}

Fix $\eps > 0$, and let $r>0$ be a constant to be specified later. By Proposition \ref{band} there is an $M = M(\eps)$ such that for all $m \geq M$, 
\beq 
    \frac{\#\{\lambda \vdash m : v(\mf f_\lambda) < r + \log_\ell m\}}{p(m)} < \eps.
\eeq
Let $M_0$ be another constant to be specified later, and put
\beq
M'=t(t+1)\max \left( M, M_0 \right).
\eeq
Let $k = tu$ and $N = kq^{k+1}M'$. Let $n \geq N$ and  $\mathsf F \in \ms F_n(\Phi)$.

Suppose $d(\mathsf F) \leq k$. Since there are $q^{k+1}$ polynomials of degree at most $k$ in $\F_q[x]$, we have $|\supp \mathsf F| \leq q^{k+1}$. Thus

\beq
    \begin{split}
        \frac{1}{(t+1)|\supp \mathsf F|}\sum\limits_{f,i} \mathsf F_i(f) &\geq \frac{\sum\limits_{f \in \Phi} |\mathsf F(f)|}{(t+1)q^{k+1}} \\
        &\geq \frac{\sum\limits_{f \in \Phi} d(f)|\mathsf F(f)|}{(t+1)kq^{k+1}} \\
        &= \frac{n}{(t+1)kq^{k+1}} \\
        &\geq \frac{M'}{t+1}.
    \end{split}
\eeq

Pick $f_0 \in \Phi$ and $0 \leq i_* \leq t$ with $\mathsf F_{i_*}(f_0)$ maximal, and put
\beq
m=\left \lfloor \frac{\mathsf F_{i_*}(f_0)}{t}\right \rfloor.
\eeq 
By the above we have $m \geq \frac{M'}{t(t+1)} \geq M$.

 For any $\bm \mu \mapsto \mathsf F$, the inequality \eqref{needed} gives  
\beq
 \left \lfloor \frac{|\core_t \bm \mu(f_0)|}{t} \right \rfloor+ \sum_{i=0}^{t-1} v(\mf f_{\bm \mu(f_0)^i}) \leq v(d_{\bm \mu});
\eeq
in particular each term on the left is a lower bound for $v(d_{\bm \mu})$.

First suppose that $i_* \neq t$. Then
\begin{equation} \label{blumb}
    \frac{\#\{\bm \mu: \bm \mu \mapsto \mathsf F, \text{ and }  v(d_{\bm \mu})  <  v((n/t)_u)\}}{\#\{\bm \mu: \bm \mu \mapsto \mathsf F\}} \leq      \frac{\#\{\bm \mu: \bm \mu \mapsto \mathsf F, \text{ and }  v(\mf f_{\mu(f_0)^{i_*}})  <  v((n/t)_u)\}}{\#\{\bm \mu: \bm \mu \mapsto \mathsf F\}}.
    \end{equation}
Since the fibre over $\mathsf F$ can be identified with the product
\beq
\prod_{f \in \supp \mathsf F} \left( \left( \prod_{i=0}^{t-1} \Lambda_{\mathsf F_i(f)} \right) \times \mc C^t_{\mathsf F_t(f)} \right),
\eeq
the right hand side of \eqref{blumb} equals
\begin{equation} \label{argos}
\frac{\#\{\lambda \vdash m : v(\mf f_\lambda) <  v((n/t)_u)\}}{p(m)}. 
\end{equation}
Now $v((n/t)_u) \leq \lceil u \rceil-\log_\ell t+ \log_\ell n$, and $n \leq t(t+1)kq^{k+1}m$, so if we put
\beq
r=\lceil u \rceil+ \log_\ell((t+1)kq^{k+1}),
\eeq
then \eqref{argos} is no greater than
\beq
\frac{\#\{\lambda \vdash m : v(\mf f_\lambda) <  r+ \log_\ell m \}}{p(m)} < \eps,
\eeq 
as desired.

Next suppose that $i_*=t$; then 
\beq
\begin{split}
    \frac{\#\{\bm \mu: \bm \mu \mapsto \mathsf F, \text{ and }  v(d_{\bm \mu})  <  v((n/t)_u)\}}{\#\{\bm \mu: \bm \mu \mapsto \mathsf F\}}  
 &=\frac{\# \left\{\lambda \in \mc C_{ \mathsf F_t(f_0)}^t : m=\left \lfloor  \frac{|\la|}{t} \right \rfloor<  r + \log_\ell m \right\}}{\#\mc C_{ \mathsf F_t(f_0)}^t}. \\ 
 \end{split}
\eeq
Now pick $M_0$ so that  $m>r+ \log_\ell m$ whenever $m \geq M_0$. Then the above numerator is $0$.
 
Next, suppose that $d(\mathsf F)>k$, say $f_1 \in \supp  \mathsf F$ with $d(f_1)>k$. We will show that for all ${\bm \mu} \mapsto  \mathsf F$ we have
\begin{equation} \label{onyx}
v(d_{\bm \mu}) \geq v((n/t)_u),
\end{equation}
 so that the numerator is again $0$. We have
\beq
\begin{split}
\frac{n}{t} &= \sum_f d(f) \frac{|{\bm \mu}(f)|}{t} \\
		&=\frac{d(f_1) |\mathsf F(f_1)|}{t} + \sum_{f \neq f_1} \frac{d(f) |\mathsf F(f)|}{t} \\
		& \geq \frac{d(f_1)-h_{f_1}|\mathsf F(f_1)|}{t} + \sum_{f  } \frac{h_f |\mathsf F(f)|}{t}, \\
		\end{split}
		\eeq
so that
\beq
\lfloor n/t \rfloor -\sum_f \left \lfloor\frac{h_f|\mathsf F(f)|}{t} \right \rfloor \geq \frac{d(f_1)-h_{f_1}}{t} -1  > \lceil u \rceil.
\eeq
It follows that
\beq
 \frac{\lfloor n/t \rfloor!}{\prod\limits_{f\in \Phi}\left \lfloor\frac{| \mathsf F(f)|h_f}{t} \right \rfloor!}  
\eeq
is an integer multiple of $(n/t)_u$, whence \eqref{onyx} holds in view of \eqref{needed}. \end{proof}

\begin{prop}\label{bottle3}
     For all $u,\eps > 0$, there is an $N = N(\eps)$ such that for all $n \geq N$,
    \beq
        \frac{\#\{\bm \mu \in \mf X_{n}: v(d_{\bm \mu}) <  v((n/t)_{u})\}}{\#\mf X_{n}} < \eps.
    \eeq
\end{prop}

\begin{proof}

 For any $\eps>0$, let $N=N(\eps)$ be as in Lemma  \ref{blanket2}. Then
    \beq
    \begin{split}
        \frac{\#\{\bm \mu \in \mf X_{n}: v(d_{\bm \mu}) <  v((n/t)_{u})\}}{|\mf X_{n}|} &= \sum_{\mathsf F } \frac{\#\{\bm \mu \mapsto \mathsf F: v(d_{\bm \mu}) <  v((n/t)_{u})\}}{|\mf X_{n}|} \\ 
       &< \eps\sum_{\mathsf F }\frac{\#\{ \bm \mu:\bm \mu \mapsto \mathsf F\}}{|\mf X_{n}|} \\
        &= \eps,
    \end{split}
    \eeq
    the sums being over $\mathsf F \in \mc I_n$.

\end{proof}

 \subsection{$\ell$-Divisibility of Character Values}
 
 Recall that $\ell$ is an odd prime   coprime to $q$, and we have fixed $g \in G_{n_0}$.

\begin{prop} \label{beep} Let $r$ be a positive integer. Then
\beq 
        \lim\limits_{n \to \infty} \frac{\#\{\bm \mu \in \mf X_{n}:  \ell^r \nmid \chi_{\bm \mu}(g)\}}{|\mf X_{n}|} = 0.
    \eeq
\end{prop}

\begin{proof}
	Now LTE gives us
    \beq
        v\left(\prod\limits_{i=0}^{n_0-1} q^{n-i}-1\right) \leq  \left \lceil \frac{n_0}{t}   \right \rceil \tau +v ((n/t)_{n_0/t}).
    \eeq
    Again Corollary \ref{helmet} implies that if $\ell^r \nmid \chi_{\bm \mu}(g)$, then 
    \beq
        v(d_{\bm \mu}) < r + v\left(\prod\limits_{i=0}^{n_0-1} q^{n-i}-1\right)  \leq  r + \left \lceil \frac{n_0}{t}   \right \rceil \tau+ v((n/t)_{n_0/t}).
    \eeq
    If we put $r_0 = r +  \left \lceil \frac{n_0}{t}   \right \rceil  \tau$, then
    \beq
    v(d_{\bm \mu}) < r_0 + v((n/t)_{n_0/t}) \leq v((n/t)_{n_0/t+\ell r_0})
    \eeq
     whenever $\ell^r \nmid \chi_{\bm \mu}(g)$. Now we deduce the conclusion from Proposition \ref{bottle3}.
\end{proof}


It is clear that Theorem \ref{nuclear} follows from Propositions \ref{bip} and \ref{beep}.

\section{A Case of Low Divisibility}  \label{def.char.sec} 

\begin{lemma}\label{desk} 
For all $n$, we have $|\mf X_{n}|\leq q^n$.
\end{lemma}

\begin{proof}
    This follows from \cite[Proposition~3.5]{fulman}, since  $\mf X_{n}$ is in bijection with the set of conjugacy classes in $G_n$.
\end{proof}

\begin{proof}[Proof of Theorem \ref{voltas}]
	If $d$ and $q = p^e$ are not coprime, then $p$ divides $d$. As a result, we have
\beq
	\frac{\#\{\bm \mu \in \mf X_n: d \mid d_{\bm \mu} \}}{|\mf X_n|}  \leq \frac{\#\{\bm \mu \in \mf X_n: p \mid d_{\bm \mu} \}}{|\mf X_n|}.
\eeq

	Let $\bm \mu \in \mf X_{n}$. From \eqref{remote} we  compute 
    \beq
        v_p(d_{\bm \mu}) = e \sum\limits_{f \in \Phi} d(f) \alpha(\bm \mu(f)).
    \eeq
    In particular, $v_p(d_{\bm \mu}) = 0$ exactly when $\alpha(\bm \mu(f)) = 0$ for all $f \in \Phi$, or equivalently, when $\bm \mu(f) = (|\bm \mu(f)|)$ for all $f \in \Phi$.  (Meaning each $\bm \mu(f)$ is a partition with one part.) Hence the correspondence
    \beq
    \bm \mu \longleftrightarrow \prod \limits_{f \in \Phi} f^{|\bm \mu(f)|}
    \eeq
    is a bijection between the set of $\bm \mu$ for which $d_{\bm \mu}$ is not divisible by $p$ and the set of monic polynomials of degree $n$  not divisible by $x$. There are $q^n - q^{n-1}$ such polynomials, so
\beq
 \frac{\#\{\bm \mu \in \mf X_n: d \mid d_{\bm \mu} \}}{|\mf X_n|} \leq \frac{|\mf X_n|-(q^n-q^{n-1})}{|\mf X_n|}  \leq \frac{1}{q}.
  \eeq 
\end{proof}

\begin{remark} Using  \cite[Proposition~3.5]{fulman}, one can similarly show
\beq
\lim_{n \to \infty}  \frac{\#\{\chi \in \irr(G_n) : p \mid \deg \chi \}}{|\irr(G_n)|} =\frac{1}{q}.
\eeq
\end{remark}

\bibliographystyle{alpha}
\bibliography{ref.bib}
\end{document}